\documentclass[12pt]{amsart}
\usepackage{amsmath,amssymb}
\usepackage{graphicx,color}
\usepackage{amscd,amsthm,amsfonts,boxedminipage}
\usepackage{boxedminipage}
\usepackage{array}

\input xy
\xyoption{all}

\newcommand{\R}{{\mathbb{R}}}
\newcommand{\Z}{{\mathbb{Z}}}
\newcommand{\C}{{\mathbb{C}}}

\newcommand{\F}{{\mathbb{F}}}

\newcommand{\bea}{\begin{eqnarray}}
\newcommand{\eea}{\end{eqnarray}}

\newcommand{\bp}{\begin{pmatrix}}
\newcommand{\ep}{\end{pmatrix}}

\newcommand{\bps}{\begin{smallmatrix}}
\newcommand{\eps}{\end{smallmatrix}}

\newcommand{\ti}{\tilde}

\def \cE{{\mathcal E}}

\def \cI{{\mathcal I}}

\def \cO{{\mathcal O}}

\def \cV{{\mathcal V}}

\def \GT{\mathcal{GT}}

\def \m{{\frak m}}

\def \0{{\bf 0}}
\def \1{{\bf 1}}
\def \ii{{\bf i}}

\def \DG{\mathit{DG}}

\def \fpartial#1{\frac{\partial}{\partial {#1}}}
\def \fpart#1#2{\frac{\partial #1}{\partial #2}}

\def \({\left(}
\def \){\right)}

\def \ti#1{\tilde{#1}}

\def \cO{{\mathcal O}}

\def \grad{\mathrm{grad}}

\def \Int{\mathrm{Int}}

\def \Tri{\mathit{Tr}}

\def \Mo{\mathit{Mo}}

\def \GT{\mathcal{GT}}
\def \HGT{\mathcal{HGT}}

 \newtheorem{thm}{Theorem}[section]
 \newtheorem{lem}[thm]{Lemma}
 \newtheorem{prop}[thm]{Proposition}
 \newtheorem{cor}[thm]{Corollary}

\theoremstyle{definition}

\begin{document}

\title{Homological mirror symmetry of $\F_1$ via Morse homotopy}

\date{\today}

\noindent{
 }\\

\author{Masahiro Futaki}
\author{Hiroshige Kajiura}
\address{Graduate School of Science, Chiba University, 
  263-8522 Japan}
\email{futaki@math.s.chiba-u.ac.jp}
\email{kajiura@math.s.chiba-u.ac.jp}

\begin{abstract}
This is a sequel to our paper 
\cite{fut-kaj1}, where 
we proposed a definition of the Morse homotopy 
of the moment polytope of toric manifolds.
Using this as the substitute of the Fukaya category of the 
toric manifolds, 
we proved a version of homological mirror symmetry 
for the projective spaces and their products 
via Strominger-Yau-Zaslow construction 
of the mirror dual Landau-Ginzburg model. 

In this paper we go this way further 
and extend our previous result to the case of 
the Hirzebruch surface $\F_1$.
\end{abstract}

\renewcommand{\thefootnote}{}
\footnote{
M.~F.\ is supported by Grant-in-Aid for Scientific Research (C) (18K03269)  
of the Japan Society for the Promotion of Science. 
H.~K.\ is supported by Grant-in-Aid for Scientific Research (C) (18K03293)  
of the Japan Society for the Promotion of Science. }

\renewcommand{\thefootnote}{\arabic{footnote}}
\setcounter{footnote}{0}

\maketitle

\tableofcontents

 \section{Introduction}

In \cite{SYZ}, Strominger-Yau-Zaslow proposed a construction 
of mirror dual Calabi-Yau manifolds via dual torus fibrations 
on a closed manifold.
Kontsevich-Soibelman \cite{KoSo:torus} proposed a framework to systematically prove 
homological mirror symmetry 
by interpolating a variant of the Morse homotopy.
Morse homotopy was first introduced for closed manifolds as a Morse theoretic 
(dimensionally reduced) model of the Fukaya category by Fukaya \cite{fukayaAinfty}.
He and Oh showed that Morse homotopy 
fully faithfully embeds into the Fukaya category of the cotangent bundles \cite{FO}.
It is hard to extend 
Kontsevich-Soibelman's program 
to general Calabi-Yau's because of the existence of singular fibers 
of the SYZ fibration.
See for example Fukaya \cite{fukaya:asymptotic} 
for an outline of the whole program on this line.

SYZ picture is also applicable to the case of 
toric Fano's and their Landau-Ginzburg mirrors, 
in which the ends of the total space of the Landau-Ginzburg model 
corresponds to the toric divisors.
It was 
first discussed in mathematical-symplectic geometric context by Auroux \cite{auroux}.
Abouzaid \cite{abouzaid09} formulated and proved a version of 
the homological mirror symmetry for toric Fano's 
by employing the tropical geometric setting 
on the Fukaya-Morse side, and 
proved an $A_\infty$ embedding by using 
an abstract framework of \v{C}ech category 
on the complex side.

Based on the differential geometric formulation by 
Leung-Yau-Zaslow \cite{LYZ} and Leung \cite{leung05}, 
Fang \cite{fang08} studied homological mirror symmetry 
for $\C P^n$ using the mirror transform associated with the SYZ fibration.
Chan \cite{chan09} studied the case of more general projective torics to determine 
which Lagrangians in the Landau-Ginzburg mirror 
correspond to holomorphic line bundles. 

We further investigate this kind of formulation to more directly realize 
the SYZ picture in the toric Fano cases in our previous paper \cite{fut-kaj1}.
We proposed a definition of the Morse homotopy $\Mo(P)$ 
for the moment polytope $P$ 
and proved a version of the homological mirror symmetry 
for projective spaces and their products.
This enables us to give a concrete description 
of the functorial mirror transform.

We go this way further in this paper and compute the case of 
the Hirzebruch surface $\F_1$. 
Namely, we prove a version of the homological mirror symmetry:
\begin{thm}[Corollary \ref{cor:prod2}]
We have an equivalence of triangulated categories
\begin{equation*}
Tr(\Mo_\cE (P) ) \simeq D^b (coh (\F_1)) ,
\end{equation*}
where $P$ is the moment polytope of $\F_1$ and 
$\cE$ is the collection of Lagrangian sections 
mirror to the chosen full strongly exceptional collection of 
holomorphic line bundles on $\F_1$.
\end{thm}

This paper is organized as follows.
In Section 2, we recall some basic settings and definitions 
from our previous paper \cite{fut-kaj1} without going into details.
In Section \ref{ssec:F_k}, we recall Hirzebruch surfaces in homogeneous coordinates.
In Section \ref{ssec:lb-F_k}, we realize holomorphic line bundles on $\F_1$ 
in a geometric way and construct the DG category $DG(\F_1)$. 
The corresponding Lagrangian sections are obtained explicitly in 
subsection \ref{ssec:Lag-ab}. 
In subsection \ref{ssec:Mo-F_1}, we compute the Morse homotopy and prove the main theorem.

While we do some calculations for $\F_k$ with general $k \geq 1$, 
we restrict ourselves in this paper to the case $k=1$ 
since $\F_k$ is Fano if and only if $k=0$ or $1$, where 
$\F_0=\C P^1\times\C P^1$ is already discussed in \cite{fut-kaj1}. 
In fact, in our setting, we mainly keep in mind that the toric manifold is Fano, 
where any line bundle is guaranteed to be an exceptional object 
by the Kodaira vanishing theorem. 
However, we should note that 
a similar equivalence may exist for $\F_k$ with $k\ge 2$ 
since it is known that 
$D^b(coh(\F_k))$ has full strongly exceptional collections of line bundles 
\cite{hille-perling11}. 
Another reason for us to restrict the case $k = 1$ in this paper is that 
we can construct the category $\Mo(P)$ and enjoy the equivalence above 
explicitly in this case. 
See subsection \ref{ssec:Lag-ab} for more details.  

\noindent
{\bf Acknowledgements} We would like to thank Fumihiko Sanda 
for telling us some known facts about 
the full strongly exceptional collections on Hirzebruch surfaces.

\section{Preliminaries on the SYZ fibrations and homological mirror symmetry}
\label{sec:SYZ}

In this section, we review some notions and settings 
from our previous paper \cite{fut-kaj1}.

\subsection{Hessian manifold and dual torus fibrations}

Let $B$ be a tropical affine manifold: 
it is equipped with an affine open cover $\{ U_\lambda \}_{\lambda \in \Lambda}$ 
whose transition functions are affine with integral linear part.
We assume for simplicity that all nonempty 
intersections of $U_\lambda$'s are contractible.
Namely, 
the coordinate transformation is of the form
\begin{equation*}
x_{(\mu)} = \varphi_{\lambda \mu} x_{(\lambda)} + \psi_{\lambda \mu}
\end{equation*}
with $\varphi_{\lambda \mu} \in GL(n;\Z)$ and $\psi_{\lambda \mu} \in \R^n$, 
where $x_{(\lambda)} = (x^1_{(\lambda)},\dots,x^n_{(\lambda)})^t$ 
and $x_{(\mu)} = (x^1_{(\mu)},\dots,x^n_{(\mu})^t$ denote 
the local coordinates on $U_\lambda$ and $U_\mu$ respectively.
We omit the suffix $(\lambda)$ when no confusion may occur.

We call $B$ Hessian when it is equipped with a metric $g$ 
locally expressed as 
\begin{equation*}
g_{ij} = \dfrac{\partial^2 \phi}{\partial x^i \partial x^j}
\end{equation*}
for some smooth local function $\phi$.
Hereafter we assume that $B$ is Hessian.

Using the metric $g$ we first define the dual affine coordinates 
on the base space as follows: 
since $\sum_{j=1}^n g_{ij} dx^j$ is closed if $(B,g)$ is Hessian, 
there exists a function $x_i := \phi_i$ of $x$ for each $i$ such that 
\begin{equation}
dx_i = \sum_{j=1}^n g_{ij} dx^j . \label{dual-coord}
\end{equation}
We thus obtain the dual coordinates 
$\check{x}^{(\lambda)} := (x_1^{(\lambda)},...,x_n^{(\lambda)})^t$.

We then denote the fiber coordinates on 
$T^\ast U_\lambda = T^\ast B|_{U_\lambda}$ dual to $x_{(\lambda)}$ 
by $(y^{(\lambda)}_1,...,y^{(\lambda)}_n)$. 
Denote by $(y_{(\lambda)}^1,...,y_{(\lambda)}^n)$ 
the fiber coordinates on $TB|_{U_\lambda}$ which corresponds to 
$(y^{(\lambda)}_1,...,y^{(\lambda)}_n)$ via the isomorphism 
$TB \stackrel{\cong}\to T^\ast B$ induced by $g$.
The cotangent bundle 
$T^\ast B$ is equipped with the standard symplectic form 
$\omega_{T^\ast B} := \sum_{i=1}^n dx^i \wedge dy_i$. 
The tangent bundle 
$TB$ is a complex manifold, where $z^i = x^i + \ii y^i$'s form the complex coordinates. 
We can further 
equipped $TB$ with the symplectic form 
$\omega_{TB} := \sum_{i,j=1}^n g_{ij} dx^i \wedge dy^j$ 
and $T^*B$ with the complex structure given by the complex coordinates $z_i = x_i + \ii y_i$'s. 
These structures in turn give the 
K\"ahler structures on both $TB$ and $T^\ast B$. 

Next we consider $\Z^n$-actions on $TB$ and $T^\ast B$.
The action of $(0,,,0.\overset{i}{1},0,...,0) \in \Z^n$ is defined by 
$y^i \mapsto y^i+2\pi$ and $y_i \mapsto y_i+2\pi$ respectively.
This is well-defined because $B$ is affine and 
the linear part $\varphi_{\lambda\mu}$ of the transition functions are integral.
Therefore we can divide $TB$ and $T^\ast B$ by this action of $\Z^n$ 
to get a pair of K\"ahler manifolds $M = TB/\Z^n$ and $\check{M} = T^\ast B/\Z^n$, 
and dual torus fibrations:
\[
\xymatrix{
M \ar[dr]_{\pi} & & \check{M} \ar[dl]^{\check{\pi}} \\
& B & . 
}
\]

\subsection{The categories $DG$ and $\cV$}

Now let $X$ be a smooth compact toric manifold.
Consider the complement of toric divisors 
$\check{M} := X \setminus \mu^{-1}(\partial P)$ for the moment map 
$\mu \colon X \to P \subset \R^n$.
By fixing an appropriate structure of the Hessian manifold 
on $B := \Int P$, we get an affine torus fibration on 
$\check{M} \to B$ whose K\"ahler structure coincides with 
that coming from the given one on $X$.
Applying the construction in the previous subsection, 
we get the dual torus fibration $M \to B$ 
with a K\"ahler structure on the total space.
We thus get the structure of affine torus fibration $\check{M} \to B$ 
and its dual torus fibration $M \to B$ with K\"ahler structures.
We describe such affine structure of $\F_k$ concretely in Section \ref{ssec:F_k}.

We consider the correspondence between Lagrangian submanifolds in $M$ 
and holomorphic vector bundles with $U(1)$-connection on $\check{M}$.
The following calculation is based on a version of the Fourier-Mukai transform 
and came from Leung-Yau-Zaslow \cite{LYZ} and Leung \cite{leung05} 
(but with different conventions).
Let $L$ be a Lagrangian section 
$\underline{y} \colon B \to M$.
Then $\underline{y}$ can 
locally expressed as $df$ for some smooth function $f$.
We associate to it a line bundle $V$ on $\check{M}$ 
with the $U(1)$-connection 
\begin{equation}
D := d - \dfrac{\ii}{2\pi} \sum_{i=1}^n y^i(x)dy_i . \label{D}
\end{equation}
This is holomorphic 
since $\underline{y}$ is a Lagrangian section.

We first define the category $\cV$ associated with $\check{M}$.
It is the DG category consisting of pairs of holomorphic line bundles and 
$U(1)$-connections of the form (\ref{D}).
More precisely, for objects $y_a = (V_a,D_a)$ and $y_b = (V_b,D_b)$ 
it has the hom space
\begin{equation*}
\cV(y_a,y_b) := \Gamma (V_a,V_b) \underset{C^\infty(\check{M})}{\otimes} 
\Omega^{0,\ast}(\check{M})
\end{equation*}
where $\Gamma (V_a,V_b)$ denotes the space of homomorphisms from $V_a$ to $V_b$.
This space is $\Z$-graded with the degree of the anti holomorphic differential forms 
and we denote the degree $r$ part by $\cV^r(y_a,y_b)$.
Decompose $D_a = D_a^{(1,0)}+D_a^{(0,1)}$ and set $d_a := 2D_a^{(0,1)}$.
The differential $d$ on $\cV(y_a,y_b)$ is then defined as
\begin{equation*}
d_{ab}(\psi) := d_b\psi - (-1)^r\psi d_a
\end{equation*}
for $\psi \in \cV^r(y_a,y_b)$.
The product structure is given by combining the composition of 
bundle homomorphisms and the wedge product:
\begin{equation*}
m(\psi_{ab},\psi_{bc}) := (-1)^{r_{ab},r_{bc}}\psi_{bc}\wedge\psi_{ab} .
\end{equation*}

We then define the DG category $DG(X)$ of holomorphic line bundles 
on the toric manifold $X$.
For a line bundle $V$ on $X$, we take a holomorphic connection $D$ 
whose restriction to $\check{M}$ is isomorphic to
a line bundle on $\check{M}$ with a connection of the form
\[
 d-\dfrac{\ii}{2\pi} \sum_{i=1}^n y^i(x) dy_i . 
\]
We set the objects of $\DG(X)$ as such pairs $(V,D)$. 
The space $\DG(X)(y_a,y_b)$ of morphisms is 
defined as the graded vector space 
each graded piece of which is given by 
\[
 \DG^r(X)(y_a,y_b):=\Gamma(V_a,V_b)\otimes\Omega^{0,r}(X) 
\]
with $\Gamma(V_a,V_b)$ being the space of smooth 
bundle morphism from $V_a$ to $V_b$. 
The composition of morphisms is defined in a similar way as
that in $\cV(\check{M})$ above.
The differential 
\[
 d_{ab}: \DG^r(X)(y_a,y_b)\to \DG^{r+1}(X)(y_a,y_b)
\]
is defined by
\[
 d_{ab}(\ti\psi)
 := 2\left(D_b^{0,1}\ti\psi- (-1)^r\ti\psi D_a^{0,1}\right) . 
\]

We then have a faithful embedding 
$\cI \colon DG(X) \to \cV$ by restricting line bundles on $X$ to $\check{M}$.
We define $\cV'$ to be the image $\cI (DG(X))$ of $DG(X)$ under $\cI$.

For a full exceptional collection $\cE$ of $DG(X)$, 
we denote the corresponding full subcategories consisting of $\cE$ 
by $DG_\cE (X) \subset DG(X)$ 
and $\cV'_\cE \subset \cV'$ respectively.

\subsection{The Morse homotopy $\Mo(P)$}

In the case of the moment polytopes, the objects of $\Mo(P)$ are 
Lagrangian sections $\underline{y} \colon B \to M$ 
which corresponds to objects of $DG(X)$ described in the previous subsection.
We shall see explicitly later in Section \ref{ssec:Lag-ab} in the case of $\F_1$.
Note that (i) they intersect cleanly, i.e.\ 
there exists an open set $\widetilde{B}$ such that $\bar{B}\subset\widetilde{B}$
and $L, L'$ over $B$ can be extended 
to graphs of smooth sections 
over $\widetilde{B}$ so that they intersect cleanly, and 
(ii) for each $L$, we can locally take a Morse function $f_L$ on $\widetilde{B}$
so that $L$ is the graph of $df_L$. 

For a given pair $(L,L')$,
we assign a grading $|V|$ for each connected component $V$
of the intersection $\pi(L\cap L')$ in $P=\bar{B}$
as the dimension of the stable manifold $S_v\subset \widetilde{B}$
of the gradient vector field $-\grad(f_L-f_{L'})$ with a point $v\in V$. 
This does not depend on the choice of the point $v\in V$. 
The space $\Mo(P)(L,L')$ of morphisms 
is then set to be the $\Z$-graded vector space
spanned by the connected components $V$ of $\pi(L\cap L')\in P$ 
{\em such that there exists a point $v\in V$ which is an interior point 
of $S_v\cap P$}. 
\footnote{We consider the Morse cohomology degree instead of the Morse homology degree.}
Note that, by this definition, the space $\Mo(P)(L,L)$ is generated by 
$P$, which is of degree zero and forms the identity morphism for any object $L$.

Rather than going into full details, 
we only explain $\m_2$ because of the following reasons: 
firstly, the Morse homotopy for $\F_1$ is minimal, i.e.\ with zero differential.
Secondly, the set of objects $\cE$ we compute later forms 
strongly exceptional collection in $Tr(\Mo_\cE (P))$ 
and therefore we do not need to compute $\m_k$ with $k\geq 3$ 
to compute $Tr(\Mo_\cE (P))$.
(For more comments, see \cite{fut-kaj1}, Section 4.5.)

Take a triple $(L_1,L_2,L_3)$, 
connected components of the intersections 
$V_{12} \subseteq L_1 \cap L_2, \ V_{23} \subseteq L_2 \cap L_3, \ 
V_{13} \subseteq L_1 \cap L_3$ and define 
$\GT(v_{12},v_{23};v_{13})$ to be the set of the trivalent gradient trees 
starting at $v_{12} \in V_{12}, \ v_{23} \in V_{23}$ and ending at $v_{13} \in V_{13}$.
Define 
$\GT (V_{12},V_{23};V_{13}) := \cup_{v_{12} \in V_{12}, v_{23} \in V_{23}, 
v_{13} \in V_{13}} \GT(v_{12},v_{23};v_{13})$ and 
$\HGT(V_{12},V_{23};V_{13}) := \GT (V_{12},V_{23};V_{13})/{\rm smooth\ homotopy}$.
This set becomes a finite set when $|V_{13}| = |V_{12}|+|V_{23}|$ and therefore we 
define $\m_2 = $ composition of morphisms by
\begin{align*}
\m_2 \colon \Mo(P)(L_1,L_2) \otimes \Mo(P)(L_2,L_3) \to \Mo(P)(L_1,L_3) \\
(V_{12},V_{13}) \mapsto \sum_{|V_{13}| = |V_{12}|+|V_{23}|} 
\sum_{[\gamma ]\in\HGT(V_{12},V_{23};V_{13})} e^{-A(\gamma )}V_{13} 
\end{align*}
where $A(\gamma )$ is the symplectic area of disk
obtained by lifting the gradient tree $\gamma$ to $M$.

 \section{Homological mirror symmetry of $\F_1$}
\label{F_1}

Following
Hille-Perling\cite{hille-perling11},
Elagin-Lunts\cite{elagin-lunts15},
Kuznesov\cite{kuz17}, etc, 
the triangulated category
$D^b(coh(\F_1))\simeq \Tri(DG(\F_1))$ 
has a series of full strongly exceptional collections 
\[
 \cE:=(\cO, \cO(1,0),\cO(c,1),\cO(1+c,1) ) , 
\]
where $\cO(a,b)$ is a line bundle on $\F_1$ we shall define later
in subsection \ref{ssec:lb-F_k}. 

We denote the corresponding full subcategories by
$DG_\cE(\F_1)\subset DG(\F_1)$, 
$\cV'_\cE\subset \cV'=\cV'(\F_1)$ and
$\Mo_\cE(P)\subset\Mo(P)$, where
$P$ is the moment polytope of $\F_1$. 

Then our main theorem is stated as follows.
\begin{thm}\label{thm:main}
For any fixed $c=0,1,\dots$, 
there exists a linear $A_\infty$-equivalence 
\[
 \iota:\Mo_\cE(P)\to \cV'_\cE
\]
such that
for any generator $V\in\Mo_\cE(P)(L,L')$ with any $L,L'\in\Mo_\cE(P)$ 
\begin{itemize}
 \item $\iota(V)\in (\cV')^0(\iota(L),\iota(L'))\subset C^\infty(B)$
extends to a continuous function on $P=\bar{B}$ and 
 \item we have 
\[
 \max_{x\in P} |\iota(V)(x)|=1,\quad \{x\in P\ |\ |\iota(V)(x)|=1 \}=V .
\]
\end{itemize}
\end{thm}
As in the cases for $X=\C P^n$ and $X=\C P^m\times\C P^n$, 
this theorem implies a version of homological mirror symmetry of $\F_1$. 
\begin{cor}\label{cor:prod1}
We have a linear $A_\infty$-equivalence
\[
 \Mo_\cE(P)\simeq \DG_\cE(\F_1) . 
\]
\qed\end{cor}
\begin{cor}\label{cor:prod2}
We have an equivalence of triangulated categories 
\[
 \Tri(\Mo_\cE(P))\simeq D^b(coh(\F_1)) 
\]
where $\Tri$ denotes the twisted complexes construction 
by Bondal-Kapranov \cite{BK:enhanced} and Kontsevich \cite{kon94}.
\qed\end{cor}

As a biproduct of the proof of Theorem \ref{thm:main},
we also show the following.
\begin{prop}\label{prop:main}
If $L\ne L'$, any generator $V\in\Mo_{\cE}(P)(L,L')$
belongs to the boundary $\partial(P)$. 

For given bases
$V\in\Mo_{\cE}(L,L')$ and
$V'\in\Mo_{\cE}(L',L'')$, the image $\gamma(T)$ by 
any gradient tree
$\gamma\in\GT(V,V';V'')$ with $V''\in\Mo_\cE(P)(L,L'')$ 
belongs to the boundary $\partial(P)$
unless $L=L'=L''$. 
\end{prop}

In subsection \ref{ssec:F_k},
we explain how to treat Hirzebruch surfaces in our set-up. 
In subsection \ref{ssec:lb-F_k}, 
we discuss line bundles on $\F_k$ constructed from
the toric divisors and construct 
the DG category $DG(\F_k)$ consisting of these line bundles. 
In subsection \ref{ssec:DG-F_1}, we calculate
the cohomology of the DG-category $DG(\F_1)$ of line bundles 
and in particular full subcategories $DG_\cE(\F_1)$ consisting of
full strongly exceptional collections $\cE$ of
the triangulated category $\Tri(DG(\F_1))\simeq D^b(coh(\F_1))$ 
following the known technique in toric geometry. 
In subsection \ref{ssec:Lag-ab}, 
we construct the Lagrangin sections which are SYZ mirror dual
to the line bundles in $DG(\F_1)$
based on a geometric realization of the line bundles
in subsection \ref{ssec:lb-F_k}. 
The obtained Lagrangian sections will be the objects of $\Mo(P)$. 
In subsection \ref{ssec:cV'-F_1}, 
we translate the cohomology $H(DG_\cE(\F_1))$ to
the cohomology $H(\cV'_\cE)$. 
In subsection \ref{ssec:Mo-F_1}, we construct $\Mo_\cE(P)$ 
and show our main theorem by comparing the result with $H(\cV'_\cE)$. 
At the end, we discuss a morphism of degree one in
$\cV'$ and see that we have the corresponding morphism in $\Mo(P)$
in subsection \ref{ssec:exmp-degree1}. 

 \subsection{Hirzebruch surfaces}
\label{ssec:F_k}

Though we need $\F_1$ only, in this subsection we treat $\F_k$
with general $k \geq 1$ since the SYZ mirror of $\F_k$
is obtained for any $k$ in a similar way.

The Hirzebruch surface $\F_k$ is defined by
\[
 \F_k:=\{[s_0:s_1], [t_0:t_1:t_2]\ |\ (s_0)^kt_0=(s_1)^k t_1\} \subset \C P^1\times \C P^2 . 
\]
If $a\ne 0$, then $t_0=(s_1/s_0)^k t_1$.
Thus, for $a\ne 0$ and $x\ne 0$, 
we have a chart 
\[
 U_1:=\{[1:u], [u^k:1:v]\}
\]
where $u=s_1/s_0$ and $v=t_2/t_1$. 
Similarly, for $s_0\ne 0$ and $t_2\ne 0$, 
we have a chart 
\[
 U_3:=\{[1:u], [u^kv:v:1]\}
\]
where $u=s_1/s_0$ and $v=t_1/t_2$. 
Similarly, 
we have a chart
\[
 U_2:=\{[u:1], [1:u^k:v]\},\quad u=s_0/s_1,\ \ v=t_2/t_0
\]
for $s_1\ne 0$ and $t_0\ne 0$, and a chart
\[
 U_4:=\{[u:1],[v:u^kv:1]\},\quad u=s_0/s_1,\ \ v=t_0/t_2
\]
for $s_1\ne 0$ and $t_2\ne 0$. 
Thus, one sees that $\F_k\simeq P(\cO(-k)\oplus \cO)$, 
where $\cO(-k)$ and $\cO$ are the line bundles over 
$\C P^1 = \{ [s_0:s_1] \}$. 

We have the natural projections
\begin{equation*}
  \xymatrix{
    & \F_k \ar[ld]_{\pi_1}\ar[rd]^{\pi_2} & \\
   \C P^1 & & \C P^2 \ \ ,
  } 
\end{equation*}
where $\pi_1$ gives the fibration structure of $P(\cO(-k)\oplus \cO)$ 
over $\C P^1$.

A K\"ahler form $\omega$ is then obtained by
\[
 \omega = C_1\pi_1^*(\omega_{\C P^1}) + C_2\pi_2^*(\omega_{\C P^2}), 
\]
where $C_1>0$ and $C_2>0$ are real constants and
$\omega_{\C P^n}$ is the Fubini-Study form
(which we treat explicitly in \cite[section 2.2]{fut-kaj1}). 
Correspondingly,
we have the moment map $\mu:\F_k\to\R^2$ defined by
\[
\mu([s_0:s_1],[t_0:t_1:t_2])
:=2\left(C_1\frac{|s_0|^2}{|s_0|^2+|s_1|^2}+C_2 k\frac{|t_1|^2}{|t_0|^2+|t_1|^2+|t_2|^2},
C_2\frac{|t_2|^2}{|t_0|^2+|t_1|^2+|t_2|^2},\right) . 
\]
The image $\mu(\F_k)$ is the trapezoid surrounded by
the $x^1$-axis, $x^2$-axis, $x^2=2C_2$ and $x^1+kx^2=2(C_2k+C_1)$. 
Namely, the moment polytope is
\[
 P:=\{ (x^1,x^2)\in\R^2\ |\ 0\le x^1\le 2(C_2k+C_1)-kx^2,\ 0\le x^2\le 2C_2\} . 
\]

Now, we treat 
\[
 \check{M}:=U_1\cap U_2\cap U_3\cap U_4 \to \Int P =: B
\]
as a torus fibration. We express this with the coordinates of $U_2=:U$,
where $\mu((u=0,v=0))=(0,0)$. 
When we further denote $u=e^{x_1+\ii y_1}$ and $v=e^{x_2+\ii y_2}$, 
then $(y_1,y_2)$ is the coordinates of a fiber of $\check{M}$. 
We see that the restriction of $\mu$ to $\check{M}\subset\F_k$
gives the fibration structure 
$\check{\pi} := \mu|_{\check{M}} \colon \check{M} \to B = \Int P$.

We further see that the dual coordinates of $(x_1,x_2)$ is
$(x^1,x^2)$ above. 
Actually, the K\"ahler form $\omega$ is expressed as 
\[
\omega = -2\ii d\left(C_1\frac{\bar{u}du}{1+u\bar{u}}
 +C_2\frac{\bar{u}^kdu^k+\bar{v}dv}{1+(u\bar{u})^k+v\bar{v}}\right) ,
\]
so one has
\[
 g^{-1}=
 4\bp
 C_1\frac{s}{(1+s)^2}+C_2\frac{k^2\cdot s^k(1+t)}{(1+s^k+t)^2} &
 C_2\frac{-k\cdot s^kt}{(1+s^k+t)^2} \\
 C_2\frac{-k\cdot s^kt}{(1+s^k+t)^2} &
 C_2\frac{t(1+s^k)}{(1+s^k+t)^2} \\
 \ep , 
\]
where $s:=u\bar{u}\ (=e^{2x_1})$, $t:=v\bar{v}\ (=e^{2x_2})$. 
By (\ref{dual-coord}), the dual coordinates are defined by 
\[
\bp dx^1\\ dx^2\ep
 = g^{-1}
  \bp dx_1 \\ dx_2\ep 
= 4\bp
 C_1\frac{s}{(1+s)^2}+C_2\frac{k^2\cdot s^k(1+t)}{(1+s^k+t)^2} &
 C_2\frac{-k\cdot s^kt}{(1+s^k+t)^2} \\
 C_2\frac{-k\cdot s^kt}{(1+s^k+t)^2} &
 C_2\frac{t(1+s^k)}{(1+s^k+t)^2} \\
 \ep
 \bp dx_1 \\ dx_2\ep, 
\]
which is satisfied by 
\begin{equation}\label{dual-coordx1x2}
 \begin{split}
   (x^1,x^2)&= \left(C_1\frac{2e^{2x_1}}{1+e^{2x_1}}
   +C_2k\frac{2 e^{2kx_1}}{1+e^{2kx_1}+e^{2x_2}},
   C_2\frac{2 e^{2x_2}}{1+e^{2kx_1}+e^{2x_2}} \right)  \\
   & = \mu([e^{x_1+\ii y_1}:1],[1:e^{k(x_1+\ii y_1)}:e^{x_2+\ii y_2}]) .  \\
 \end{split}
\end{equation}
Hereafter we fix $C_1=C_2=1$ since 
the structure of the category $\Mo(P)$ we shall construct
is independent of these constants.

 \subsection{Line bundles on $\F_k$}
\label{ssec:lb-F_k}

Any line bundle over $\F_k$ is constructed from 
a toric divisor, 
which is a linear combination of the following four divisors 
\begin{equation*}
 \begin{split}
   & D_{12}=(t_2=0) \\
   & D_{24}=(s_0=t_1=0) \\
   & D_{13}=(s_1=t_0=0) \\
   & D_{34}=(t_0=t_1=0) . 
 \end{split}
 \end{equation*}
Note that $D_{24}$ is the fiber of $\pi_1$ at $[s_0:s_1] = [0:1]$, 
$D_{13}$ is the fiber of $\pi_1$ at $[1:0]$, 
and the remaining two are sections of $\pi_1$.
The corresponding Cartier divisors are as follows.
\begin{equation*}
 \begin{split}
   D_{12}: & \{(U_1,v_1),(U_2,v_2),(U_3,1),(U_4,1)\}\\
   D_{24}: & \{(U_1,1),(U_2,u_2),(U_3,1),(U_4,u_4)\} \\
   D_{13}: & \{(U_1,u_1),(U_2,1),(U_3,u_3),(U_4,1)\} \\
   D_{34}: & \{(U_1,1),(U_2,1),(U_3,v_3),(U_4,v_4)\} . 
 \end{split}
 \end{equation*}
Now, the coordinate transformations are
\begin{equation*}
  \begin{split}
 (u_1,v_1)& =(1/u_2,v_2/u_2^k), \\
 (u_3,v_3)&=(1/u_2,u_2^k/v_2), \\
 (u_4,v_4)&=(u_2,1/v_2) .   
  \end{split}
\end{equation*}
The transition functions are then
\begin{equation*}
 \begin{array}{cccc}
   D_{12}: & \phi_{21}=v_1/v_2=u_2^{-k},& \phi_{23}= 1/v_2,& \phi_{24}=1/v_2 , \\
   D_{24}: & \phi_{21}=1/u_2,& \phi_{23}=1/u_2,& \phi_{24}=u_4/u_2=1, \\
   D_{13}: & \phi_{21}=u_1=1/u_2,& \phi_{23}=u_3=1/u_2,& \phi_{24}=1, \\
   D_{34}: & \phi_{21}=1,& \phi_{23}=v_3=u_2^k/v_2,& \phi_{24}=v_4=1/v_2  .   
 \end{array}
 \end{equation*}
Thus, we see that $D_{13}$ and $D_{24}$ define the same line bundle; $\cO(D_{13})=\cO(D_{24})$, 
and $\cO(D_{12}) = \cO(D_{34} + k D_{24})$.  
In this sense, 
any line bundle over $\F_k$ 
is generated either by $(D_{24},D_{34})$ or by $(D_{24},D_{12})$.

On the other hand, 
we have line bundles $\pi_1^*(\cO_{\C P^1}(1))$ 
and $\pi_2^*(\cO_{\C P^2}(1))$ over $\F_k$ 
via 
the pair of projections
\begin{equation*}
  \xymatrix{
    & \F_k \ar[ld]_{\pi_1}\ar[rd]^{\pi_2} & \\
   \C P^1 & & \C P^2 \ \ .   
  } 
\end{equation*}
The transition functions for them
are obtained by the pullbacks of 
$\cO_{\C P^1}(1)\to\C P^1$ and $\cO_{\C P^2}(1)\to\C P^2$ by
$\pi_1^\ast$ and $\pi_2^*$, respectively. 
Then we can identify 
$\pi_1^*(\cO_{\C P^1}(1))=\cO(D_{24})$ and $\pi_2^*(\cO_{\C P^2}(1))=\cO(D_{12})$.

The connection one-forms for 
$\pi_1^*(\cO_{\C P^1}(1))$ and $\pi_2^*(\cO_{\C P^2}(1))$
are also obtained by the pullbacks. 
On $U_2$, they are expressed as  
\begin{equation*}
 \begin{split}
   \pi_1^*\left(A_{\C P^1}\right)
   & = -\frac{\bar{u}du}{1+u\bar{u}}  , \\
   \pi_2^*\left(A_{\C P^2}\right)
   & = -\frac{\bar{u}^k d(u^k)+\bar{v}dv}{1+(u\bar{u})^k+v\bar{v}} , 
 \end{split}
\end{equation*}
where $A_{\C P^n}$ is the connection one-form
for the line bundle $\cO_{\C P^n}(1)$ 
on $\C P^n$ as is presented explicitly
for instance in \cite[section 3.3]{fut-kaj1}. 
We again denote $s:=u\bar{u}$, $t:=v\bar{v}$, and then 
\begin{equation*}
 \begin{split}
   \pi_1^*\left(A_{\C P^1}\right)
   & = -\frac{s(dx_1+\ii dy_1)}{1+s}  , \\
   \pi_2^*\left(A_{\C P^2}\right)
   & = -\frac{k s^k(dx_1+\ii dy_1)+t(dx_2+\ii dy_2)}{1+s^k+t} .
 \end{split}
\end{equation*}
Similarly, the connection one-forms for
$\cO(a,b):=\cO(a D_{24}+b D_{12})$ is given by
\[
 A_{(a,b)}:= -a\frac{s(dx_1+\ii dy_1)}{1+s} -b \frac{k s^k(dx_1+\ii dy_1)+t(dx_2+\ii dy_2)}{1+s^k+t} .
\]
In a similar way as for $\C P^n$ in \cite[section 3.3]{fut-kaj1}, 
the $dx_1$ term and the $dx_2$ term are removed by the isomorphisms as
follows.
\begin{equation}\label{Psi-ab}
 \begin{split}
 &\Psi_{(a,b)}^{-1} (d + A_{(a,b)})\Psi_{(a,b)} = d 
-a\ii\frac{s dy_1}{1+s} -b\ii\frac{k s^k dy_1 +t dy_2}{1+s^k+t} , \\ 
 & \Psi_{(a,b)}
 := \left(1+s\right)^{\frac{a}{2}}\left(1+s^k+t\right)^{\frac{b}{2}}  .
 \end{split}
\end{equation}
Note that the connection in the right hand side above
is of the form (\ref{D}). 
Thus, we set $DG(\F_k)$
as the DG-category consisting of the line bundles $\cO(a,b)$
with any $(a,b)$.

 \subsection{Cohomologies of the DG category $DG(\F_1)$}
\label{ssec:DG-F_1}

We first discuss the global sections of line bundles $\cO(a,b)$ over $\F_k$. 
The structure of the global sections is obtained by 
following \cite[p.66]{fulton93toric}. 
In particular, we see that
\[
 d_{(a,b)}:=\dim\left(\Gamma(\F_k,\cO(a,b))\right) = 
  (a+1) + (a+1+k) + \cdots + (a+1+kb) \\
  = \frac{(b+1)(2a+2+kb)}{2} .  
\]
For instance, $d_{(1,0)}=2$, $d_{(0,1)}=2+k$, $d_{(a,0)}=a+1$.  
Using the coordinates $(u,v)$ for $U_2$, the generators
of $\Gamma(\cO(a,b))$ are expressed as 
\begin{equation}\label{DGgen}
 \begin{array}{ccccc}
  u^0 , & u^1, & \dots, & \dots, & u^{a+kb}, \\
  u^0v^1, & u^1v^1, & \dots, & u^{a+k(b-1)}v^1 ,&  \\
   \vdots, & \vdots &   & & \\
  u^0v^b, & \dots & u^av^b & &  .   
 \end{array}
\end{equation}
Namely, 
\[
 \psi_{(i_1,i_2)}:= u^{i_1}v^{i_2}
\]
are the generators, where $0\le i_2\le b$ and $0\le i_1\le a+k(b-i_2)$. 

Since each $\cO(a,b)$ is a line bundle, 
we have
\[
 DG^0(\F_k)(\cO(a_1,b_1),\cO(a_2,b_2))
 \simeq DG^0(\F_k)(\cO,\cO(a_2-a_1,b_2-b_1)) . 
\]
Thus, we obtain any zero-th cohomology of the space of morphisms in $DG(\F_k)$ 
from $\Gamma(\cO(a,b))$, $a,b\in\Z$, i.e.,
\[
H^0(DG(\F_k)(\cO(a_1,b_1),\cO(a_2,b_2)))\simeq
H^0(DG(\F_k)(\cO,\cO(a_2-a_1,b_2-b_1)))\simeq 
\Gamma(\F_k, \cO(a_2-a_1,b_2-a_1)) .
\]

We can also calculate the cohomologies 
\[
H^r(DG(\F_k)(\cO(a_1,b_1),\cO(a_2,b_2)) \simeq
H^r(DG(\F_k)(\cO,\cO(a_2-a_1,b_2-b_1))
\]
with $r>0$ in the way written in \cite[p.74]{fulton93toric}. 
Now we consider $\F_1$. It is known
(Hille-Perling \cite{hille-perling11}) 
that 
\[
 \cE:=(\cO, \cO(1,0),\cO(c,1),\cO(1+c,1) )
\]
with a fixed $c=0,1,\dots$
forms a full strongly exceptional collection of
$\Tri(DG(\F_1))\simeq D^b(coh(\F_1))$. 
Let us consider the full subcategory $DG_\cE(\F_1)\subset DG(\F_1)$ 
consisting of $\cE$. 
We already calculated the zero-th cohomologies of the space of morphisms 
in $DG_\cE(\F_1)$. 
We can also check that 
we have no nontrivial cohomologies of the space of morphisms of degree $r>0$
in $DG_\cE(\F_1)$. 
These calculations give a direct confirmation of the fact that 
$\cE$ is actually a strongly exceptional collection in $\Tri(DG_\cE(\F_1))$, 
and agree, for instance, 
with the Euler bilinear form on $K_0(D^b(coh(\F_1)))$ 
presented in Kuznesov \cite[Example 3.7]{kuz17}.

 \subsection{Lagrangian sections $L(a,b)$}
\label{ssec:Lag-ab}

We continue to concentrate on the case $\F_1$ and 
let us discuss the Lagrangian section $L(a,b)$
in the fiber $M\to B$ corresponding to the line bundle $\cO(a,b)$.  
Comparing the connection one-form in (\ref{Psi-ab})
with (\ref{D}), we see that $L(a,b)$ is expressed as
the graph of
\[
\bp y^1 \\ y^2 \ep = 2\pi
\bp
a\frac{s}{1+s}+ b\frac{s}{1+s+t} \\ 
b\frac{t}{1+s+t} 
\ep
\]
where $s=e^{2x_1}$ and $t=e^{2x_2}$. 
Now, let us rewrite $s,t$ by $x^1,x^2$. Recall that
the dual coordinates are given in (\ref{dual-coordx1x2}):
\begin{align}
  x^1&=\frac{2s}{1+s}+\frac{2s}{1+s+t} \label{x1} , \\
  x^2&=\frac{2t}{1+s+t} . \label{x2}
\end{align}
By (\ref{x2}), $t$ is expressed as 
\[
 t=\frac{x^2(1+s)}{2-x^2} . 
\]
Substituting this to (\ref{x1}) yields
\[
x^1 =\frac{2s}{1+s}+ \frac{2s}{(1+s)+\frac{x^2}{2-x^2}(1+s)}
    = \frac{(4-x^2)s}{1+s}
\]
and hence we obtain
\[
 \frac{s}{1+s}= \frac{x^1}{4-x^2} . 
\]
So, $s=x^1/(4-x^1-x^2)$, $1+s= (4-x^2)/(4-x^1-x^2)$,
and we get $t=(1+s)x^2/(2-x^2)= (4-x^2)x^2/(2-x^2)(4-x^1-x^2)$. 
To summarize, we obtain 
\begin{equation*}
 \begin{split}
 \frac{s}{1+s} &= \frac{x^1}{4-x^2}, \\
 \frac{s}{1+s+t} &= \frac{x^1(2-x^2)}{2(4-x^2)}, \\
 \frac{t}{1+s+t} &= \frac{x^2}{2} . 
 \end{split}  
\end{equation*}
In particular, the Lagrangian $L(a,b)$ corresponding to
$\cO(a,b):=\cO(a D_{24}+b D_{12})$ is given by
\begin{equation*}
\bp y^1 \\ y^2 \ep
=
\bp a\frac{x^1}{4-x^2}+b\frac{x^1(2-x^2)}{2(4-x^2)} \\
    b\frac{x^2}{2}
\ep
=
\bp
 \frac{(2a+(2-x^2)b) x^1}{2(4-x^2)} \\ \frac{bx^2}{2}
\ep
. 
\end{equation*}
The corresponding Morse function $f$ is given by
\begin{equation}\label{f-F1}
 \begin{split}
 f & = \frac{a}{2}\log(1+s) +\frac{b}{2}\log(1+s^k+t) \\
 & = +\frac{1}{2}\log
 \left(\frac{4-x^2}{4-x^1-x^2}\right)^a
 \left(\frac{2(4-x^2)}{(2-x^2)(4-x^1-x^2)}\right)^b . 
 \end{split}
\end{equation}
For $\F_k$ with general $k>1$, the Lagrangian section $L(a,b)$
should still be obtained in a similar way. 
However, we do not seem to obtain a closed formula for $(s, t)$ in terms of $(x^1,x^2)$.

 \subsection{Cohomologies $H(\cV'_\cE)$}
\label{ssec:cV'-F_1}

We consider $DG(\F_1)$ and the the faithful embedding
$\cI: DG(\F_1)\to\cV$ where $\cV=\cV(\check{M})$ is the DG category
of line bundles on $\check{M}$. The image is denoted by 
$\cV':=\cI(DG(\F_1))$. 
Then, 
each generator $\psi_{(i_1,i_2)}$ in (\ref{DGgen}) is sent to be 
\begin{equation}\label{pre-e}
\Psi_{(a,b)}^{-1} \psi_{(i_1,i_2)} =
 (4-x^1-x^2)^{\frac{a+b-i_1-i_2}{2}}(2-x^2)^{\frac{b-i_2}{2}}
(4-x^2)^{-\frac{a+b-i_2}{2}}(x^1)^{\frac{i_1}{2}}(x^2)^{\frac{i_2}{2}}
 e^{\ii (i_1 y_1+i_2y_2)} 
\end{equation}
in $\cV'$. 
Namely, these form a basis of
$H^0(\cV'(\cO,\cO(a,b)))$. 
A basis of $H^0(\cV'(\cO(a_1,b_1),\cO(a_1+a),\cO(b_1+b)))$
is of the same form.
We rescale each basis $\Psi_{(a,b)}^{-1} \psi_{(i_1,i_2)}$ by multiplying
a positive number and denote it by ${\bf e}_{(a,b);(i_1,i_2)}$
so that 
\[
 \max_{x\in P} |{\bf e}_{(a,b);(i_1,i_2)}(x)| = 1 . 
\]
Note that ${\bf e}_{(a,b);(i_1,i_2)}$
is a function on $B$, but can be extended continuously
to that on $P=\bar{B}$
since the exponents  
\[
 \frac{a+b-i_1-i_2}{2},\quad
 \frac{b-i_2}{2},\quad
 \frac{i_1}{2},\quad
 \frac{i_2}{2}
\]
in (\ref{pre-e})are non-negative. 
This shows the former statement about the properties of $\iota(V)$
in Theorem \ref{thm:main}.

 \subsection{$\Mo_\cE(P)$}
\label{ssec:Mo-F_1}

The objects of $\Mo(P)$ are the Lagrangian sections $L(a,b)$ 
obtained in subsection \ref{ssec:Lag-ab}. 
Since we have
\[
 \Mo(P)(L(a_1,b_1),L(a_2,b_2))
 \simeq \Mo(P)(L(0,0),L(a_2-a_1,b_2-b_1)), 
\]
we concentrate on calculating the space $\Mo(P)(L(0,0),L(a,b))$. 
We discuss that when there exists a nonempty intersection of
\[
\bp y^1\\y^2\ep
=
\bp y^1_{(a,b)}(x)\\y^2_{(a,b)}(x)\ep
=
\bp
 \frac{(2a+(2-x^2)b) x^1}{2(4-x^2)} \\ \frac{bx^2}{2}
 \ep
\]
with 
\[
 \bp y^1\\y^2\ep = \bp i_1\\ i_2\ep 
\]
in the covering space of $\bar{M}\to P$. 

We first consider $L(a,b)$ with $b\ge 0$. 
Since $0\le x^2\le 2$, we have
\[
   0\le \frac{bx^2}{2}=i_2\le b . 
\]
If we further assume $a+b-i_2\ge 0$, then 
we also have
\begin{equation}\label{cd-y_1}
0\le i_1=\frac{(2a+(2-x^2)b)x^1}{2(4-x^2)} \le
 \frac{(2a+(2-x^2)b)}{2}= a+b-i_2 ,  
\end{equation}
where we used $4-x^2\ge x^1$ in the inequality. 
We will discuss the case $a+b-i_2<0$ later. 
By solving $y^j_{(a,b)}(x)=i_j$, $j=1,2$, we obtain the following. 
\begin{lem}\label{lem:pre-gen-Mo}
We assume $b\ge 0$ and $(a,b)\ne (0,0)$. 
For any $(i_1,i_2)$ satisfying
\[
 0\le i_2\le b,\quad 0\le i_1\le a+b-i_2, 
\]
the intersection $V_I$ is nonempty and connected. 
\begin{itemize}
  \item If $b\ne 0$ and $a+b-i_2\ne 0$, then $V_I$ consists of the point $v_I$ 
such that 
\[
 x^1(v_I)=\frac{4-2i_2/b}{a+b-i_2}i_1, \qquad
 x^2(v_I)=\frac{2i_2}{b} .
\]
 \item If $b=0$, then 
$i_2=0$ and then the intersection is 
\[
 V_{(i_1,0)}:=\{(x^1,x^2)\in P\ |\ x^1=\frac{4-x^2}{a} i_1\} .
\]
 \item If $a+b-i_2=0$, then $i_1=0$ and 
the intersection is
\[
 V_{(0,a+b)}:= \{(x^1,2 + \frac{2a}{b})\in P\} .
\]
\end{itemize}
\qed\end{lem}

Note that the condition $a+b-i_2=0$ is satisfied only when $a\le 0$ and $b>0$.
\begin{lem}\label{lem:gen-Mo}
We assume $b\ge 0$ and $(a,b)\ne (0,0)$.   
For any $I=(i_1,i_2)$ satisfying
\[
 0\le i_2\le b,\quad 0\le i_1\le a+b-i_2, 
\]
the intersection $V_I$ forms a generator of
$\Mo(P)(L(a_1,b_1),L(a_1+a,b_1+b))$ of degree zero. 
\end{lem}
\begin{pf}
The gradient vector field associated to $V_I$ is of the form
\begin{equation}\label{grad1}
 \left(\frac{(a+b-bx^2/2)x^1}{4-x^2}-i_1\right)\fpartial{x^1}
 + \left(\frac{bx^2}{2}-i_2\right)\fpartial{x^2} .
\end{equation}
If $b\ne 0$ and $a+b-i_2\ne 0$, then
$V_I$ consists of the point $v_I$, and 
the stable manifold $S_{v_I}$ of the gradient vector field
is $\{v_I\}$ itself, so $V_I$ is a generator of degree zero. 
If $b=0$, the gradient vector field associated to $V_{(i_1,0)}$ is
of the form
\begin{equation}\label{grad2}
\left(\frac{ax^1}{4-x^2}-i_1\right)\fpartial{x^1} . 
\end{equation}
For each $v\in V_{(i_1,0)}$, the stable manifold $S_v$ of the
gradient vector field is again $\{v\}$ itself,
so $V_{(i_1,0)}$ gives a generator of degree zero. 
If $a+b-i_2=0$, then the gradient vector field
associated to $V_{(0,a+b)}$ is of the form 
\begin{equation*}\label{grad3}
 (a+b-bx^2/2)\left( \frac{x^1}{4-x^2}\fpartial{x^1}
 -\fpartial{x^2}\right) .
\end{equation*}
Because of the term $-(a+b-bx^2/2)\partial/\partial x^2$, 
for each $v\in V_{(0,a+b)}$, the stable manifold $S_v$ of the
gradient vector field is again $\{v\}$ itself, 
so $V_{(0,a+b)}$ gives a generator of degree zero. 
\qed\end{pf}

Now, we discuss the case $a+b-i_2<0$.
This case occurs only if $a$ is negative
since $a<-(b-i_2)\le 0$. 
We have
\[
 a+b-i_2\le i_1\le 0
\]
instead of (\ref{cd-y_1}). 
The intersection $V_I$ again consists of 
a point $v_I$ given by 
\[
 x^2=\frac{2i_2}{b},\quad
 x^1=\frac{4-2i_2/b}{a+b-i_2}i_1
\]
if $b\ne 0$, and 
\[
 V_{(i_1,0)}:=\{(x^1,x^2)\in P\ |\ x^1=\frac{4-x^2}{a} i_1\} 
\]
if $b=0$. 
\begin{lem}\label{lem:non-gen}
We assume $b\ge 0$ and $(a,b)\ne (0,0)$.   
For any $I=(i_1,i_2)$ satisfying
\[
 0\le i_2\le b,\quad a+b-i_2\le i_1\le 0, 
\]
we consider the intersection $V_I$ which is nonempty and connected.
\begin{itemize}
 \item[(i)] The intersection does not form 
a generator of
$\Mo(P)(L(a_1,b_1),L(a_1+a,b_1+b))$ of degree zero. 
 \item[(ii)] The intersection does not form 
a generator of
$\Mo(P)(L(a_1,b_1),L(a_1+a,b_1+b))$ (o any degree) if 
$a= -1$.  
\end{itemize}
\end{lem}
\begin{pf}
The gradient vector field associated to $V_I$ is also of the same form
as in the proof of Lemma \ref{lem:gen-Mo}.  
If $b>0$, then the gradient vector field is (\ref{grad1}) but
now the sign of the coefficient for $x^1$ is reversed compared to
the case $a+b-i_2>0$. 
Thus, we have $|V_I|=1$. 
Similarly, if $b=0$, then the gradient vector field is (\ref{grad2}), 
but again the sign of the coefficient for $x^1$ is reversed compared to
the case $a+b-i_2>0$ and we have $|V_I|=1$. 
Thus, the statement (i) is proved. 

We consider the case $a=-1$, where $b=i_2$ holds since $a+b-i_2<0$ and
$i_1=0,-1$. 
If $b>0$, then $V_I$ consists of the point
$v_I=(-2i_1,2)$. 
We see that the stable manifold is 
\[
 S_{v_I}=\{x^2=2\} .
\]
in both cases $I=(0,b)$ and $I=(-1,b)$. 
Thus, $v_{(0,b)}=(0,2)$ (resp. $v_{(-1,b)}=(2,2)$) 
is not an interior point of $S_{v_{(0,1)}}\cap P\subset S_{v_{(0,1)}}$
(resp. $S_{v_{(-1,b)}}\cap P\subset S_{v_{(-1,b)}}$). 
If $b=0$, then we have $V_{(0,0)}=D_{24}$ and $V_{(-1,0)}=D_{13}$. 
For $v\in D_{24}$, the stable manifold is 
\[
 S_v=\{x^2=x^2(v)\} , 
\]
so $v$ is not an interior point of $S_v\cap P \subset S_v$. 
Similarly, any $v\in D_{13}$ is not an interior point of
$S_v\cap P \subset S_v$. 
Thus, any $V_I$ does not form a generator of
$\Mo(P)(L(a_1,b_1),L(a_1+a,b_1+b))$ if $a=-1$. 
\qed\end{pf}

By Lemma \ref{lem:gen-Mo} and Lemma \ref{lem:non-gen} (i), 
we see that each generator $V_I$ of degree zero is in one-to-one
correspondence with the generator ${\bf e}_{(a,b);I}$ in subsection \ref{ssec:cV'-F_1}. 
\begin{lem}\label{lem:morp-corresp}
We assume $b\ge 0$ and $(a,b)\ne (0,0)$. 
Each generator ${\bf e}_{(a,b);I}\in H^0(\cV'(\cO(a_1,b_1),\cO(a_1+a,b_1+b)))$ 
is expressed as the form
\[
 {\bf e}_{(a,b);I}(x)= e^{-f_I} e^{\ii I y}, 
\]
where $f_I$ is the $C^\infty$ function on $P$ satisfying
\[
df_I= \sum_{j=1}^2\fpart{f_I}{x_j} dx_j ,\quad
 \fpart{f_I}{x_j}=y_{(a,b)}^j-i_j
\]
in $B$ and $\min_{x\in P}f_I(x)=0$. 
In particular, we have 
\[
 \{x\in P\ |\ f_I(x)=0\}= V_I .
\]
\end{lem}
\begin{pf}
The statement of the first half is guaranteed by our construction
$DG(\F_1)\simeq\cV'\subset\cV$. 
Of course, we can also check it directly. 
By (\ref{f-F1}), the function $f_I$
satisfying $\fpart{f_I}{x_j}=y_{(a,b)}^j-i_j$, $j=1,2$,
is given by 
\[
\frac{1}{2}
\log\left(\frac{4-x^2}{4-x^1-x^2}\right)^a\left(\frac{2(4-x^2)}{(2-x^2)(4-x^1-x^2)}\right)^b  -i_1 x_1 - i_2 x_2, 
\]
where $i_1,i_2 \in\Z$. 
Since we have 
\begin{equation*}
 \begin{split}
   x_1& = \frac{1}{2}\log s
   = \frac{1}{2}\log\left(\frac{x^1}{4-x^1-x^2}\right), \\
   x_2&= \frac{1}{2}\log t
   = \frac{1}{2}\log\left(\frac{x^2(4-x^2)}{(2-x^2)(4-x^1-x^2)}\right), 
 \end{split}
\end{equation*}
we obtain 
\[
f_I=-\left(\log\left((4-x^1-x^2)^{\frac{a+b-i_1-i_2}{2}}(2-x^2)^{\frac{b-i_2}{2}}
(4-x^2)^{-\frac{a+b-i_2}{2}}(x^1)^{\frac{i_1}{2}}(x^2)^{\frac{i_2}{2}}\right)
\right)
+\text{const.} . 
\]

The latter half can be shown directly by 
rewriting 
\[
 \begin{split}
  & (4-x^1-x^2)^{\frac{a+b-i_1-i_2}{2}}(2-x^2)^{\frac{b-i_2}{2}}
(4-x^2)^{-\frac{a+b-i_2}{2}}(x^1)^{\frac{i_1}{2}}(x^2)^{\frac{i_2}{2}}\\
 & = \left(1-\frac{x^1}{4-x^2}\right)^{\frac{a+b-i_1-i_2}{2}}
     \left(\frac{x^1}{4-x^2}\right)^{\frac{i_1}{2}}
     (2-x^2)^{\frac{b-i_2}{2}}(x^2)^{\frac{i_2}{2}} 
 \end{split}
\]
and regarding the result as a function in variables
$x^1/(4-x^2)$ and $x^2$. 
\qed\end{pf}

Now, we fix $c\ge 0$ and consider
the full subcategory $\Mo_\cE(P)\subset\Mo(P)$ consisting of
$\cE=(L(0,0),L(1,0),L(c,1),L(1+c,1))$. 
We call an element of $\Mo(P)(L,L')$
an {\bf ordered} (resp.\ {\bf non-ordered}) morphism
if $\cE=(...,L,...,L',...)$ (resp. $\cE=(...,L',...,L,...)$). 
We see that the space
\[
 \Mo_\cE(P)(L(a_1,b_1),L(a_1+a,b_1+b))\simeq \Mo(P)(L(0,0),L(a,b))
\]
of ordered morphisms satisfies $0\le b\le 1$. 
We call an ordered morphism with $b=0$ (resp.\ $b=1$)
a morphism {\bf of type} $b=0$ (resp. $b=1$). 
\begin{lem}\label{lem:morp}
In $\Mo_\cE(P)$, we have only morphisms of degree zero,
where each generator $V_I\ne P$ belongs to $\partial(P)$. 
Then, the correspondence 
$V_I\mapsto \iota(V_I)={\bf e}_{ab;I}$ gives 
a quasi-isomorphism
\[
 \iota:\Mo_\cE(P)(L(a_1,b_1),L(a_2,b_2))\to
 \cV'_\cE(\cO(a_1,b_1),\cO(a_2,b_2)) 
\]
of complexes. 
\end{lem}
\begin{pf}
By Lemma \ref{lem:non-gen} (ii),
we see that restricting the correspondence in Lemma \ref{lem:morp-corresp}
to $\Mo_\cE(P)$ yields the quasi-isomorphism of this Lemma for
each space of ordered morphisms. 

In particular, an ordered morphism of type $b=0$ has $a=1$,
and an ordered morphism of type $b=1$ satisfies $a\ge -1$. 
In both cases, 
we can check directly 
that the generators given in Lemma \ref{lem:pre-gen-Mo} 
belongs to $\partial(P)$. 

It remains to calculate non-ordered morphisms. 
Then, we may consider the opposite case to Lemma \ref{lem:pre-gen-Mo}
in the sense that we consider the case $b< 0$ and
\[
 b\le i_2\le 0,\quad a+b-i_2\le i_1\le 0 . 
\]
We see that the corresponding intersection is the same as
the one in Lemma \ref{lem:pre-gen-Mo}. Namely, we have
\[
 V_{(-a,-b);-I} = V_{(a,b);I} (=:V_I) 
\]
for $b\ge 0$. 
The sign of the corresponding gradient vector field is reversed
compared to that in the proof in \ref{lem:gen-Mo}. 
Thus, we have 
\begin{itemize}
 \item $|V_{(-a,-b);-I}|=2$ if $b>0$ and $a+b-i_2>0$,  
 \item $|V_{(-a,-b);-I}|=1$ if $b=0$,  
 \item $|V_{(-a,-b);-I}|=1$ if $a+b-i_2=0$. 
\end{itemize}
In particular, for $b=0$ or $a+b-i_2=0$, the stable manifold $S_v$ of
a point $v\in V_{(-a,-b);-I}$ intersect transversally with $V_{(-a,-b;-I)}$ at
$v$. 
Since $V_{(-a,-b);-I}\in\partial(P)$ in all these three cases, 
we can conclude that $V_{(-a,-b);-I}$ cannot be a generator 
of $\Mo_\cE(L(a_1+a, b_1+b),L(a_1,b_1))$. 
\qed\end{pf}

Now we discuss the composition structure in $\Mo_\cE(P)$. 
Let us examine the gradient flows starting from a point in 
a generator $V_I\in\Mo(P)(L(0,0),L(a,b))$ with $b=0$ and $b=1$. 
If $b=0$, then the gradient vector field is of the form
\[
 \grad f = \left(\frac{ax^1}{4-x^2}-i_1\right)\fpartial{x^1}
\]
as we already saw, where we have $i_2=0$. 
Thus, the gradient trajectories starting from $v=(x^1(v),x^2(v))\in V_I$ 
are always in the line $x^2=x^2(v)$. 
If $b=1$, then the gradient vector field is of the form 
\[
\grad f = \left(\frac{(a+1-x^2/2)x^1}{4-x^2}-i_1\right)\fpartial{x^1}
 + \left(\frac{x^2}{2}-i_2\right)\fpartial{x^2} .
\]
In this case, we see that 
$V_{(i_1,i_2=0)}$ belongs to $D_{12}\subset\partial(P)$, 
(where $V_{(0,0)}=D_{12}$ if $a=-1$ and otherwise $V_{(i_1,0)}$ consists of
the point $v_{(i_1,0)}:=(4/(a+1),0)$), 
and $V_{(i_1,i_2=1)}$ belongs to $D_{34}\subset\partial(P)$. 
(Note that $V_{(0,1)}=D_{34}$ if $a=0$,
and otherwise $V_{(i_1,1)}$ consists of a point. )
This implies that the gradient trajectories
starting from $v_{(i_1,0)}$ are always in the line $x^2=0$,
and the gradient trajectories starting from a point in $V_{(i_1,1)}$ 
are always in the line $x^2=2$.

Suppose that we have a composition in $\Mo_\cE(P)$ 
\[
 V_I\cdot W_J 
\]
of two generators such that $V_I\ne P$ and $W_J\ne P$. 
We have such a composition only when 
\begin{itemize}
 \item $V_I$ is of type $b=0$ and $W_J$ is of type $b=1$ or 
 \item $V_I$ is of type $b=1$ and $W_J$ is of type $b=0$. 
\end{itemize}
In both cases, the result $V_I\cdot W_J$ is generated by a generator $Z_{I+J}$
of type $b=1$ with index $I+J$ 
since indices are preserved by the composition.
This means that
$Z_{I+J}$ belongs to
$D_{12}$ (resp. $D_{34}$)
if the generator, $V_I$ or $W_J$, of type $b=1$
belongs to $D_{12}$ (resp. $D_{34}$).
Since the gradient trajectories starting from a point in the generator,
$V_I$ or $W_J$, of type $b=0$ run horizontally,
we see that the image $\gamma(T)$ of the gradient tree $\gamma$ defining
the product $V_I\cdot W_J$ belongs to $D_{12}$ or $D_{34}$.
This gives the proof of Proposition \ref{prop:main}.

Let $v:=V_I\cap\gamma(T)$, $w:=W_J\cap\gamma(T)$, and 
$z:=Z_{I+J}\cap\gamma(T)$. 
Note that $v$ and $w$ are
the images of the two external vertices of the trivalent tree $T$
by $\gamma$, whereas $z$ is the image of the root vertex of $T$ by $\gamma$. 
Then, $z$ sits on the interval $vw$ and
the image of the root edge of $T$ by $\gamma$ is $\{z\}$. 
If we express $\iota(V_I)=e^{-f_v}\cdot e^{\ii I y}$, 
then $f_v(v)=0$ and $f_v(z)$ is the symplectic area
of the triangle disk whose edges are the interval $vz$ and
the corresponding two Lagrangian sections on $vz$. 
Similarly, for $\iota(W_J)=e^{-f_w}\cdot e^{\ii J\check{y}}$, 
the value $f_w(z)$ is the symplectic area of the
corresponding triangle disk. 
This shows the compatibility
\[
 \iota(V_I\cdot W_J)=\iota_I(V_I)\cdot\iota(W_J) .
\]
This completes the proof of Theorem \ref{thm:main}.

 \subsection{Example of morphisms of degree one: $\Mo(L(0,0),L(2,-2))$}
\label{ssec:exmp-degree1}
 
So far, we do not see any morphism of higher degree. 
In this subsection, we calculate the space $\Mo(L(0,0),L(2,-2))$ 
and show that it includes a generator of degree one. 
The intersections of $\pi(L(\cO))$ with $\pi(L(\cO(2,-2)))$ are 
obtained by solving 
\begin{equation*}
  \begin{split}
     y^1&=-2\frac{x^1(2-x^2)}{2(4-x^2)}+2\frac{x^1}{4-x^2}
     = \frac{x^1x^2}{4-x^2} = i_1\\ 
     y^2&=-2\frac{x^2}{2}=-x^2 =i_2.
  \end{split}  
\end{equation*} 
We have the following intersections in $\pi^{-1}(P)$ as follows. 
\begin{itemize}
 \item If $y^2=i_2=0$, then $x^2=0$, where $y^1=i_1=0$ for any $x^1$. 
 \item If $y^2=i_2=-1$, then $x^2=1$, where $y^1=x^1/3=i_1$.
In this case, we have two choices $(i_1=0,x^1=0)$ and $(i_1=1,x^1=3)$. 
 \item If $y^2=i_2=-2$, then $x^2=2$, where $y^1=x^1=i_1$. 
In this case, we have three choices $i_1=x^1=0,1,2$. 
\end{itemize}
To summarize, we have
\[
 \begin{split}
 & V_{(0,0)}=\{(x^1,0) \in P\} \\
 & V_{(0,-1)}=\{(0,1)\in P\},\quad V_{(1,-1)}=\{(3,1) \in P\} \\ 
   & V_{(0,-2)}=\{(0,2)\in P\},\quad V_{(1,-2)}=\{(1,2) \in P\},
   \quad V_{(2,-2)}=\{(2,2) \in P\} .\\ 
 \end{split}
\]
We concentrate on the examples $V_{(0,-1)},V_{(1,-1)}$, which turn out
to be generators of morphisms of degree one. 
For $V_{(0,-1)}$, the associated gradient vector field is
\[
\frac{x^1x^2}{4-x^2}\fpartial{x^1} - (x^2-1)\fpartial{x^2} , 
\]
so the stable manifold turns out to be
\[
 S_{(0,1)}=\{(0,x^2)\} . 
\]
The point $(0,1)$ is an interior point of $S_{(0,1)}\cap P\subset S_{(0,1)}$,
so $V_{(0,-1)}$ turns out to be a generator of degree one. 
Similarly, the gradient vector field associated to $V_{(1,-1)}$ is
\[
 \begin{split}   &
 \left(\frac{x^1x^2}{4-x^2}-1\right)\fpartial{x^1}
   - (x^2-1)\fpartial{x^2} \\
   & = \left(\frac{(x^1-3)(x^2-1)+(x^1-3)+4(x^2-1)}{3-(x^2-1)}\right)\fpartial{x^1}
      - (x^2-1)\fpartial{x^2},  \\
 \end{split}    
\]
so we see that $S_{(3,1)}$ is of one-dimension around $(3,1)$. 
In particular, we can check that
\[
 S_{(3,1)}= \{ (x^1,x^2)\ |\ 4-x^1-x^2 =0 \}
\]
since
\[
 \left(\frac{x^1x^2}{4-x^2}-1\right)\fpartial{x^1}
 - (x^2-1)\fpartial{x^2}
 = (x^2-1)\left( \fpartial{x^1}-\fpartial{x^2}\right) 
\]
on $4-x^1-x^2=0$. 
Thus, $V_{(1,-1)}$ also forms a generator of degree one. 

These results of course agree with the structure of $DG(\F_1)$ where $H^1(DG(\F_1)(\cO(0,0),\cO(2,-2)))$ is of two dimension.

\end{document}